\documentclass{article}

\usepackage[T1]{fontenc}

\usepackage{graphicx}
\usepackage{amsmath}
\usepackage{amsfonts}
\usepackage{amssymb}
\usepackage{amsthm}
\usepackage[all]{xy}
\input diagxy

\theoremstyle{plain}

\newtheorem*{thm}{Quillen's Theorem}
\newtheorem{theorem}{Theorem}[section]
\newtheorem{lemma}[theorem]{Lemma}
\newtheorem{proposition}[theorem]{Proposition}

\theoremstyle{definition}

\newtheorem{definition}[theorem]{Definition}
\newtheorem{notation}[theorem]{Notation}

\newtheorem{remark}[theorem]{Remark}

\def\To{\mathrel{\xy \ar@{=>} (200,0) \endxy}}
\def\comp{\raisebox{0.2mm}{\ensuremath{\scriptstyle\circ}}}

\newcommand{\fref}[1]{Figure~\ref{#1}}

\newcommand{\del}{\partial}
\newcommand{\im}{\ensuremath{\mathsf{Im\,}}}
\newcommand{\inj}{\ensuremath{\mathrm{in}}}

\renewcommand{\ker}{\ensuremath{\mathsf{Ker\,}}}
\DeclareMathOperator{\Coeq}{Coeq}

\newcommand{\coeq}{\ensuremath{\mathsf{Coeq\,}}}

\newcommand{\Ac}{\ensuremath{\mathcal{A}}}

\newcommand{\Sc}{\ensuremath{\mathcal{S}}}

\newcommand{\Z}{\ensuremath{\mathbb{Z}}}
\newcommand{\N}{\ensuremath{\mathbb{N}}}

\newcommand{\Set}{\ensuremath{\mathsf{Set}}}

\newcommand{\Gp}{\ensuremath{\mathsf{Gp}}}

\newcommand{\PCh}{\ensuremath{\mathsf{PCh}}}

\newcommand{\simp}[1]{\Sc #1}
\newcommand{\simpA}{\simp{\Ac}}

\newcommand{\cycle}{\nabla}
\newcommand{\El}{\Lambda}

\DeclareMathOperator{\pr}{pr}
\renewcommand{\hom}{\mathrm{Hom}}
\DeclareMathOperator{\op}{op}
\newcommand{\noproof}{\hfill \qed}

\newbox\pullbackbox
\setbox\pullbackbox=\hbox{\xy 0;<1mm,0mm>: \POS(4,0)\ar@{-} (0,0) \ar@{-} (4,4)
\endxy}
\def\pullback{\copy\pullbackbox}

\newbox\pushoutbox
\setbox\pushoutbox=\hbox{\xy 0;<1mm,0mm>: \POS(0,4)\ar@{-} (0,0) \ar@{-} (4,4)
\endxy}

\begin{document}

\newdir{>>}{{}*!/9pt/@{|}*!/4.5pt/:(1,-.2)@^{>}*!/4.5pt/:(1,+.2)@_{>}}
\newdir{ >>}{{}*!/9pt/@{|}*!/4.5pt/:(1,-.2)@^{>}*!/4.5pt/:(1,+.2)@_{>}}
\newdir{ >}{{}*!/-5.5pt/@{|}*!/-10pt/:(1,-.2)@^{>}*!/-10pt/:(1,+.2)@_{>}}
\newdir2{ >}{{}*-<1.8pt,0pt>@{ >}}
\newdir{^ (}{{}*!/-10pt/@{>}}
\newdir{_ (}{{}*!/-10pt/@{>}}
\newdir{>}{{}*:(1,-.2)@^{>}*:(1,+.2)@_{>}}
\newdir{<}{{}*:(1,+.2)@^{<}*:(1,-.2)@_{<}}

\title{Simplicial homotopy
in semi-abelian categories}

\author{Tim Van der Linden\footnote{Vakgroep Wiskunde, Vrije Universiteit Brussel, Pleinlaan~2, 1050~Brussel, Belgium. Email:~\texttt{tvdlinde@vub.ac.be}}}

\maketitle

\begin{abstract}
\noindent We study Quillen's model category structure for homotopy of simplicial objects in the context of Janelidze, M\'arki and Tholen's semi-abelian categories. This model structure exists as soon as $\Ac$ is regular Mal'tsev and has enough regular projectives; then the fibrations are the Kan fibrations of $\Sc\Ac$. When, moreover, $\Ac$ is semi-abelian, weak equivalences and homology isomorphisms coincide.

\smallskip
\noindent\emph{Keywords: semi-abelian category, model category, simplicial object}
\end{abstract}

\section*{Introduction}\label{Section-Introduction}

Semi-abelian categories, introduced by Janelidze, M\'arki and Tholen in~\cite{Janelidze-Marki-Tholen}, relate to the categories of groups, rings and Lie algebras in the same way as abelian categories relate to categories of modules: They form a general framework that captures many of the homological properties of the categories modelled. In any semi-abelian category, the fundamental diagram lemmas, such as the Short Five Lemma, the $3 \times 3$ Lemma, the Snake Lemma~\cite{Bourn2001} and Noether's Isomorphism Theorems~\cite{Borceux-Bourn} hold. Moreover, the usual (abelian-category) definition of homology of chain complexes still makes sense: Any short exact sequence of proper chain complexes---chain complexes of which the differentials have a normal image---induces a long exact homology sequence~\cite{EverVdL2}.

Since, via the Moore normalization functor $N\colon {\Sc \Ac \to\PCh\Ac}$, any simplicial object $A$ induces a proper chain complex $N (A)$, $A$ has homology objects $H_{n}A=H_{n}N (A)$; since $N$ is an exact functor, any short exact sequence of simplicial objects in $\Ac$ induces a long exact homology sequence in $\Ac$. These homology objects were introduced in~\cite{EverVdL2} to prove a general version of Hopf's formula---a purely homological result. On the other hand, simplicial objects are widely used for their convenient homotopical properties. The purpose of this paper is showing that, as is well-known for e.g., groups or abelian categories, homology of simplicial objects is compatible with the existing homotopy theory.

We shall give an interpretation of the following theorem due to Quillen.

\begin{thm}\cite[Theorem II.4.4]{Quillen}
Let $\Ac$ be a finitely complete category with enough regular projectives. Let $\simpA$ be the category of simplicial objects over $\Ac$. Define a map $f$ in $\simpA$ to be a fibration (resp.\ weak equivalence) if $\hom (P,f)$ is a fibration (resp.\ weak equivalence) in $\Sc \Set$ for each projective object $P$ of $\Ac$, and a cofibration if $f$ has the left lifting property with respect to the class of trivial fibrations (i.e., maps both fibration and weak equivalence). Then $\simpA$ is a closed simplicial model category if $\Ac$ satisfies one of the following extra conditions:
\begin{itemize}
\item[\emph{(*)}] every object of $\simpA$ is fibrant;
\item[\emph{(**)}] $\Ac$ is cocomplete and has a set of small regular projective generators.\noproof 
\end{itemize}
\end{thm}

The condition (**) concerns finitary varieties in the sense of universal algebra. We shall focus on categories satisfying the first condition: This class contains all abelian categories; more generally, when $\Ac$ has coequalizers of kernel pairs, (*) is satisfied exactly when $\Ac$ is Mal'tsev.

In particular, if $\Ac$ is semi-abelian with enough projectives, $\simpA$ carries a model category structure. Our aim is to describe this structure in terms of homology of simplicial objects.

\section{Fibrations are Kan fibrations}\label{Section-Fibrations}
In this short section we discuss the fibrations occurring in Quillen's model structure: In a regular Mal'tsev category with enough projectives, the fibrations are the Kan fibrations.

First of all, note that Quillen's Theorem is indeed the same as Theorem~II.4.4 in~\cite{Quillen}, where the projective objects are chosen relative to the class of so-called \emph{effective} epimorphisms. A morphism $f\colon {A\to B}$ in a finitely complete category $\Ac$ belongs to this class when for every object $T$ of $\Ac$, the diagram in $\Set$ 
\[
\xymatrix{{\hom (B,T)} \ar[rr]^-{(\cdot)\circ f} && {\hom (A,T)} \ar@<.5ex>[rr]^-{(\cdot)\circ k_{0}} \ar@<-.5ex>[rr]_-{(\cdot)\circ k_{1}} && {\hom (R[f],T)}}
\]
is an equalizer. This, however, just means that $f$ is a coequalizer of its kernel pair $k_{0},k_{1}\colon {R[f]\to A}$; hence in $\Ac$, regular epimorphisms (i.e., coequalizers) and effective epimorphisms coincide.

 From now on, we restrict ourselves to the context of \emph{regular} categories: finitely complete ones with coequalizers of kernel pairs and pullback-stable regular epimorphisms. In a regular category, the \emph{image factorization} $\im f\comp p$ of a map $f\colon {A\to B}$ as a regular epimorphism $p$ followed by a monomorphism $\im f$ is obtained as follows: $p$ is a coequalizer of a kernel pair $k_{0},k_{1}\colon {R[f]\to A}$ and the \emph{image} $\im f\colon {I[f]\to B}$ is the universally induced arrow. Now with Proposition~2 in~\cite[\S II.4]{Quillen}, Quillen shows that the effective epimorphisms in $\Ac$ are pullback-stable as soon as there are enough regular projectives, hence \emph{a finitely complete category with coequalizers of kernel pairs and enough regular projectives is regular} and asking regularity reduces to asking for the existence of coequalizers of kernel pairs. 

Recall from \cite{Carboni-Kelly-Pedicchio, EverVdL2} the following definitions. Consider a simplicial object $A$ in a regular category $\Ac$. For $n\geq 1$ and $k\in [n]$, a family
\[
x= (x_{i}\colon X\to A_{n-1})_{i\in [n],i\neq k}
\]
is called an \emph{$(n,k)$-horn} of $A$ if it satisfies $\del_{i}\comp x_{j}=\del_{j-1}\comp x_{i}$, for $i<j$ and $i,j\neq k$. A map $f \colon {A\to B}$ of simplicial objects is a \emph{Kan fibration} if, for every $(n,k)$-horn $x$ of $A$ and every $b\colon {X\to B_{n}}$ with $\del_{i}\comp b=f_{n-1}\comp x_{i}$ for all $i\neq k$, there is a regular epimorphism $p\colon {Y\to X}$ and a map $a\colon {Y\to A_{n}}$ such that $f_{n}\comp a=b\comp p$ and $\del_{i}\comp a=x_{i}\comp p$ for all $i\neq k$. A simplicial object $K$ is \emph{Kan} when the unique map ${K\to 1}$ to the terminal simplicial object $1$ is a Kan fibration.

It is easily seen that in a regular category $\Ac$ with enough regular projectives, a simplicial object $K$ is Kan if and only if for every projective object $P$ and for every $(n,k)$-horn 
\[
x= (x_{i}\colon P\to K_{n-1})
\]
of $K$, there is a morphism $y\colon P\to K_{n}$ such that $\del_{i}\comp y=x_{i}$, $i\neq k$. This means that $K$ is Kan if and only if $\hom (P,K)=\hom (P,\cdot)\comp K$ is a Kan simplicial set, for every projective object $P$. The similar statement for Kan fibrations amounts to the following:

\begin{proposition}\label{Proposition-Fibration-is-Kan-Fibration}
Let $\Ac$ be a regular category with enough projectives. A simplicial morphism $f\colon A\to B$ in $\Ac$ is a Kan fibration if and only if it is a fibration in the sense of Quillen's Theorem: For every projective object $P$ of $\Ac$, the induced morphism of simplicial sets
\[
f\comp (\cdot)=\hom (P,f)\colon \hom (P,A)\to \hom (P,B)
\]
is a Kan fibration.\noproof 
\end{proposition}

Thus we see that a regular category $\Ac$ satisfies condition (*) if and only if every simplicial object in $\Ac$ is Kan. In the introduction to~\cite{Barr}, M.~Barr conjectures that the latter condition on $\Ac$ means that every reflexive relation in $\Ac$ is an equivalence relation. Recall that this last property, together with finite completeness, gives the definition of \emph{Mal'tsev category} due to Carboni, Lambek and Pedicchio~\cite{Carboni-Lambek-Pedicchio}. In~\cite{Carboni-Kelly-Pedicchio}, Carboni, Kelly and Pedicchio prove M.~Barr's conjecture, and show that the Mal'tsev axiom is equivalent to the permutability condition $RS=SR$ for arbitrary equivalence relations $R$, $S$ on an object $X$ of~$\Ac$.

Briefly, simplicial objects in a Mal'tsev category with enough projectives and coequalizers of kernel pairs admit a model category structure, where the fibrations are just Kan fibrations. However, stronger conditions on $\Ac$ are needed for an internal description of the weak equivalences. We shall now focus on semi-abelian categories and describe the weak equivalences using homology.

\section{Homology of simplicial objects}\label{Section-Simplicial-Homology}

A \textit{semi-abelian} category is pointed, Barr exact and Bourn protomodular with binary coproducts~\cite{Janelidze-Marki-Tholen, Borceux-Bourn}. \textit{Pointed} means that it has a zero object: an initial object that is also terminal. An \textit{exact} category is regular and such that any internal equivalence relation is a kernel pair~\cite{Barr}. A pointed and regular category is \textit{protomodular} when the Short Five Lemma holds: For any commutative diagram
\[
\xymatrix{K[p'] \ar@{{ >}->}[r]^-{\ker p'} \ar[d]_-u & E' \ar@{->>}[r]^-{p'} \ar[d]^-v & B' \ar[d]^-w \\ K[p] \ar@{{ >}->}[r]_-{\ker p} & E \ar@{->>}[r]_-{p} & B}
\]
such that $p$ and $p'$ are regular epimorphisms, $u$ and $w$ being isomorphisms implies that $v$ is an isomorphism~\cite{Bourn1991}. The existence of binary coproducts entails finite cocompleteness. 

Next to the examples mentioned in the Introduction, also e.g., all abelian categories; all varieties of $\Omega$-groups (i.e., varieties of universal algebras with a unique constant and an underlying group structure), in particular the categories of commutative algebras, crossed modules and precrossed modules; and the categories of Heyting algebras, compact Hausdorff groups, non-unital $C^{*}$ algebras are well-known to be semi-abelian.

Dominique Bourn showed in~\cite{Bourn1996} that any finitely complete protomodular category is Mal'tsev; hence a semi-abelian category with enough projectives satisfies condition (*). We briefly recall from~\cite{EverVdL2}---see also~\cite{Borceux-Bourn}---some definitions and results concerning homology of simplicial objects in this context.

 From now on, $\Ac$ will always denote a semi-abelian category. A sequence of morphisms $(f_{i})_{i\in I}$
\[
\xymatrix{\dots \ar[r] & X_{i+1} \ar[r]^-{f_{i+1}} & X_{i} \ar[r]^-{f_{i}} & X_{i-1}\ar[r] & \dots}
\]
in $\Ac$ is called \emph{exact at $X_{i}$} if $\im f_{i+1}=\ker f_{i}$. It is called \emph{exact} when it is exact at $X_{i}$, for all $i\in I$. A short sequence
\begin{equation}\label{Short-Exact-Sequence}
\xymatrix{0 \ar[r] & K \ar@{{ >}->}[r]^-{k} & A \ar@{-{ >>}}[r]^-{f} & B \ar[r] & 0}
\end{equation}
is exact if and only if $f$ is a regular epimorphism (a cokernel of $k$) and $k$ is a kernel of $f$. 

A \emph{chain complex} $C$ in $\Ac$ is a sequence of morphisms $d_{n}\colon C_{n}\to C_{n-1}$, $n\in \Z$, where $d_{n}\comp d_{n+1}=0$. $C$ is called \emph{proper} when it has differentials $d_{n}$ of which the image $\im d_{n}\colon {I[d_{n}]\to C_{n-1}}$ is a kernel. As in the abelian case, the \emph{$n$-th homology object} $H_{n}C$ of a proper chain complex $C$ is the cokernel of the factorization $d'_{n+1}\colon {C_{n+1}\to K[d_n]}$ of $d_{n+1}$. Given a proper chain complex $C$, $H_{n}C$ is zero if and only if $C$ is exact at $C_{n}$; and by the Snake Lemma~\cite{Bourn2001, Borceux-Bourn}, any short exact sequence of proper chain complexes induces a long exact homology sequence.

When working with simplicial objects in $\Ac$, we shall use the notations of~\cite{Weibel}: A simplicial object $A\colon {\Delta^{\op}\to \Ac}$ consists of objects $(A_{n})_{n\in \N}$, face operators $\del_{i}\colon {A_{n}\to A_{n-1}}$ for $i\in [n]=\{0,\dots ,n \}$ and $n\in \N_{0}$, and degeneracy operators $\sigma_{i}\colon {A_{n}\to A_{n+1}}$, for $i\in [n]$ and $n\in \N$, subject to the simplicial identities. The \emph{normalization functor} $N\colon {\simpA \to \PCh \Ac}$ turns a simplicial object $A$ into the \emph{Moore complex} $N (A)$ of $A$, the chain complex with $N_{0}A=A_{0}$,
\[
N_{n} A=\bigcap_{i=0}^{n-1}K[\del_{i}\colon A_{n}\to A_{n-1}]
\]
and differentials $d_{n}=\del_{n}\comp \bigcap_{i}\ker \del_{i}\colon N_{n} A\to N_{n-1} A$, for $n\geq 1$, and $A_{n}=0$, for $n<0$. Since, by Theorem~3.6 in~\cite{EverVdL2}, $N (A)$ is a proper chain complex, one can define its homology objects as $H_{n}A=H_{n}N (A)$. The exactness of the functor $N$ \cite[Proposition 5.6]{EverVdL2} now implies:

\begin{proposition}\label{Proposition-Long-Exact-Homology-Sequence-Simp}\cite[Corollary 5.7]{EverVdL2}
Let $\Ac$ be a semi-abelian category. Any short exact sequence~\ref{Short-Exact-Sequence} in $\simpA$ canonically induces a long exact sequence
\[
\resizebox{\textwidth}{!}{\xymatrix@R=1pc{
{\dots} \ar[r] &H_{n+1}B \ar[r]^-{\delta_{n+1}} & H_nK \ar[r]^-{H_{n}k} & H_nA \ar[r]^-{H_{n}f} & H_nB \ar[r]^-{\delta_{n}} & H_{n-1}K \ar[r] & {\dots}\\
{\dots} \ar[r] &H_{1}B \ar[r]^-{\delta_{1}} & H_0K \ar[r]^-{H_{0}k} & H_0A \ar@{-{ >>}}[r]^-{H_{0}f} & H_0B \ar[r]^-{\delta_{0}} & 0}} 
\]
in $\Ac$.\noproof 
\end{proposition}

Our first target is characterizing regular epic homology isomorphisms. To do so, we need a precise description of the \emph{acyclic} simplicial objects, i.e., those $A$ in $\simpA$ satisfying $H_{n}A=0$ for all $n\in \N$.

\section{Acyclic objects}\label{Section-Acyclic-Objects}

What does it mean, for a simplicial object, to have zero homology? We give an answer to this question, not in terms of an associated chain complex, but in terms of the simplicial object itself. Recall Corollary~3.10 in~\cite{EverVdL2}:

\begin{proposition}\label{Proposition-Characterization-H_{0}}
For any simplicial object $A$, $H_{0}A=\Coeq [\del_{0},\del_{1}]$.\noproof 
\end{proposition}

\begin{figure}
\begin{center}
\resizebox{\textwidth}{!}{\includegraphics{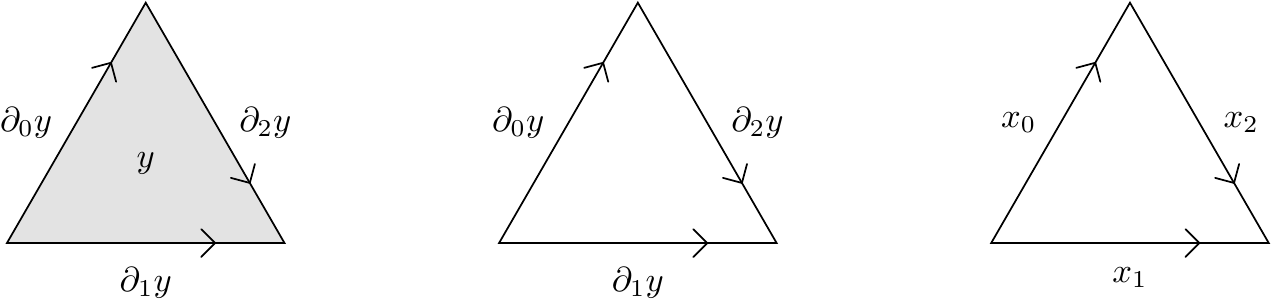}}
\end{center}
\caption{A $2$-simplex $y$, its boundary $\del y$ and a $1$-cycle $(x_{0},x_{1},x_{2})$.}\label{Figure-Acyclic-Simplicial-Object}
\end{figure}

It implies that $H_{0}A=0$ if and only if $(\del_{0},\del_{1})\colon {A_{1}\to A_{0}\times A_{0}}$ is a regular epimorphism. (Since $\Ac$ is exact Mal'tsev, the image of the reflexive graph $\del_{0},\del_{1}\colon {A_{1}\to A_{0}}$ is a kernel pair of $\coeq (\del_{0},\del_{1})\colon {A_{0}\to \Coeq [\del_{0},\del_{1}]}$; and $H_{0}A$ is zero if and only if this kernel pair is $\pr_{0},\pr_{1}\colon {A_{0}\times A_{0}\to A_{0}}$.) A~similar property holds for the higher homology objects; in order to write it down, we need the following definition (cf.\ \fref{Figure-Acyclic-Simplicial-Object}).

\begin{definition}\label{Definition-Cycles}
Let $A$ be a simplicial object in $\Ac$. The object $\cycle_{n}A$ of \emph{$n$-cycles in $A$} is defined by $\cycle_{0}A=A_{0}\times A_{0}$ and 
\[
\cycle_{n}A=\{(x_{0},\dots,x_{n+1})\,|\, \text{$\del_{i}\comp x_{j}=\del_{j-1}\comp x_{i}$, for all $i<j\in [n+1]$}\}\subset A_{n}^{n+2},
\]
for $n>0$. This gives rise to functors $\cycle_{n}\colon {\Sc \Ac \to \Ac}$. (We index projections according to the simplicial notation, i.e., $\pr_{i} (x_{0},\dots,x_{n+1})=x_{i}$ etc.)
\end{definition}

\begin{notation}\label{Notation-A^-}
For any simplicial object $A$ in $\Ac$, the simplicial object defined by $A^-_n=A_{n+1}$,  $\del^-_i=\del_{i+1}\colon {A_{n+1}\to A_n}$, and $\sigma^-_i=\sigma_{i+1}\colon {A_{n+1}\to A_{n+2}}$, for $i\in [n]$,  $n\in \N$ is denoted by $A^-$. This is the simplicial object obtained from $A$ by leaving out $A_{0}$ and, for $n\in \N$, all $\del_0\colon {A_n\to A_{n-1}}$ and $\sigma_0\colon {A_n\to A_{n+1}}$. Observe that $\del=(\del_{0})_n$ defines a simplicial morphism from $A^-$ to $A$, and write $\El A$ for its kernel. 
\end{notation}

The following result due to J.~Moore, well-known to hold in the abelian case (cf.\ Exercise 8.3.9 in~\cite{Weibel}) is easily generalized to semi-abelian categories:

\begin{lemma}\label{Lemma-Lambda}
For any $n\in \N $, $N_{n} (\El A)\cong N_{n+1}A$; hence $H_{n} (\El A)\cong H_{n+1}A$.\noproof 
\end{lemma}

We also need the following generalization of the notion of \emph{regular pushout}~\cite{Carboni-Kelly-Pedicchio}.

\begin{definition}\label{Definition-regepi-Pushout}
A square in $\Ac$ with horizontal regular epimorphisms
\[
\vcenter{\xymatrix{ A' \ar@{-{ >>}}[r]^-{f'} \ar[d]_-{v} & B' \ar[d]^-w \\
A \ar@{-{ >>}}[r]_-f & B}}
\]
is called a \emph{regular pushout} when the comparison map $(v,f')\colon  {A'\to A\times_{B}B'}$ to the pullback $A\times_{B}B'$ of $w$ along $f$ is a regular epimorphism. (The maps $v$ and $w$ are \textit{not} demanded to be regular epimorphisms.)
\end{definition}

\begin{remark}\label{Remark-Regepi-Pushout}
Every regular pushout is a pushout, but a pushout need not be a regular pushout: A counterexample may be constructed by choosing $B'=0$, $v=\inj_{A'}\colon {A'\to B+A'}$ and $f=[1_{B},0]\colon {B+A'\to B}$ in, say, the category $\Gp$ of all groups.
\end{remark}

\begin{proposition}\label{Proposition-Regular-Pushout}
In the diagram
\[
\xymatrix{0 \ar[r] & K[f'] \ar@{{ >}->}[r] \ar[d]_-u & A' \ar@{-{ >>}}[r]^-{f'} \ar[d]^-v & B' \ar[d]^-w \ar[r] & 0\\
0 \ar[r] & K[f] \ar@{{ >}->}[r] & A \ar@{-{ >>}}[r]_-{f} & B \ar[r] & 0}
\]
the right hand side square is a regular pushout if and only if the induced arrow $u$ is a regular epimorphism. 
\end{proposition}
\begin{proof}
This is a consequence of the Short Five Lemma for regular epimorphisms (\cite[Corollary 9]{Bourn1996} or \cite[Lemma~4.2.5]{Borceux-Bourn}).
\end{proof}

\begin{lemma}\label{Lemma-Cycle-and-Pullbacks}
For any $n\in \N$, the diagram
\[
\xymatrix{{\cycle_{n+1}A}\ar[d]_-{(\pr_{i+1})_{i\in [n+1]}} \ar[r]^-{\pr_{0}}& {A_{n+1}} \ar[d]^-{(\del_{i})_{i\in [n+1]}} \\
{\cycle_{n}A^{-}} \ar[r]_-{\cycle_{n}\del} & {\cycle_{n}A}}
\]
is a pullback.
\end{lemma}
\begin{proof}
Via the Yoneda Lemma, it suffices to check this in $\Set$.
\end{proof}

\begin{proposition}\label{Proposition-Characterization-Acyclic-Object}
A simplicial object $A$ in $\Ac$ is acyclic if and only if for every $n\in \N$, the morphism $(\del_{i})_{i\in [n+1]}\colon {A_{n+1}\to\cycle_{n}A}$ is a regular epimorphism.
\end{proposition}
\begin{proof}
We show by induction that, when $H_{i}A=0$ for all $i\in [n-1]$, the morphism 
\[
(\del_{i})_{i\in [n+1]}\colon {A_{n+1}\to\cycle_{n}A}
\]
is regular epic if and only if $H_{n}A=0$.

We already considered the case $n=0$ above. Hence suppose that $n\in \N$ and $H_{i}A=0$ for all $i\in [n]$, and consider the following diagram.
\[
\xymatrix{0 \ar[r] & {\El_{n+1}A} \ar@{{ >}->}[rr]^-{\ker \del_{0}} \ar[d]_-{(\del_{i}^{\Lambda})_{i\in [n+1]}} && A_{n+1}^{-} \ar[d]_-{(\del_{i+1})_{i\in [n+1]}} \ar@{-{ >>}}[rr]^-{\del_{0}} && A_{n+1} \ar[d]^-{(\del_{i})_{i\in [n+1]}} \ar[r] & 0\\
0 \ar[r] & {\cycle_{n}\El A} \ar@{{ >}->}[rr]_-{\ker \cycle_{n}\del} && \cycle_{n} A^{-} \ar[rr]_-{\cycle_{n}\del} && \cycle_{n}A \ar@{.>}[r] & 0}
\]
By the induction hypothesis on $A$, the arrow $(\del_{i})_{i}\colon {A_{n+1}\to \cycle_{n}A}$ is a regular epimorphism. It follows that so is $\cycle_{n}\del$, and both rows in this diagram are exact sequences.

Using Lemma~\ref{Lemma-Lambda} and applying the induction hypothesis to $\El A$, we get that the arrow $(\del_{i}^{\Lambda })_{i}\colon {\El_{n+1}A\to \cycle_{n}\El A}$ is a regular epimorphism if and only if $0=H_{n}\El A\cong H_{n+1}A$. Now, by Proposition~\ref{Proposition-Regular-Pushout}, we see that $H_{n+1}A=0$ exactly when the right hand side square above is a regular pushout. By Lemma~\ref{Lemma-Cycle-and-Pullbacks}, this latter property is equivalent to the induced arrow 
\[
(\del_{i})_{i\in [n+2]}\colon {A_{n+2}=A_{n+1}^{-}\to \cycle_{n+1}A}
\]
being a regular epimorphism.
\end{proof}

The meaning of this property becomes very clear when we express it as in~\fref{Figure-Acyclic-Simplicial-Object}: Up to enlargement of domain, $A$ is acyclic when every $n$-cycle is a boundary of an $(n+1)$-simplex. 

\section{Acyclic fibrations}\label{Section-Acyclic-Fibrations}

A Kan fibration is called \emph{acyclic} when it is a homology isomorphism. Here we give an explicit description of the acyclic fibrations.

Note that, using Proposition~\ref{Proposition-Long-Exact-Homology-Sequence-Simp}, regular epic homology isomorphisms may be characterized as those regular epimorphisms that have an acyclic kernel. This is why Proposition~\ref{Proposition-Characterization-Acyclic-Object} is useful in the proof of 

\begin{proposition}\label{Proposition-Characterization-Reg-Epi-Homology-Iso}
A regular epimorphism $p\colon {E\to B}$ is a homology isomorphism if and only if for every $n\in \N$, the diagram
\begin{equation}\label{Diagram-Acyclic-Fibration}
\vcenter{\xymatrix{E_{n+1} \ar[d]_-{(\del_{i})_{i\in [n+1]}} \ar@{-{ >>}}[rr]^-{p_{n+1}} && B_{n+1} \ar[d]^-{(\del_{i})_{i\in [n+1]}} \\
\cycle_{n}E \ar[rr]_-{\cycle_{n}p} && \cycle_{n}B}}
\end{equation}
is a regular pushout.
\end{proposition}
\begin{proof}
We give a proof by induction on $n$. If $n=0$ then $\cycle_{0}p$, as a product of regular epimorphisms, is regular epic. By Proposition~\ref{Proposition-Regular-Pushout}, it is now clear that Diagram~\ref{Diagram-Acyclic-Fibration} is a regular pushout if and only if $(\del_{0},\del_{1})\colon {K[p]_{1}\to \cycle_{0}K[p]}$ is a regular epimorphism, or, by Proposition~\ref{Proposition-Characterization-Acyclic-Object}, $H_{0}K[p]=0$. 

Now suppose that $n>0$. We first have to show that the arrow $\cycle_{n}p$ is a regular epimorphism. To do so, consider an arrow $x\colon {X\to \cycle_{n}B}$. By induction, a regular epimorphism $q\colon {Y\to X}$ exists and arrows $y_{k}\colon {Y\to E_{n-1}}$, such that for each element $k$ of the set
\[
K=\{\del_{j}\comp \pr_{i}\comp x\colon X\to B_{n-1}\,|\,j\in [n],\,i\in [n+1] \},
\]
$k\comp q= p_{n-1}\comp y_{k}$, and such that these $y_{k}$ assemble to arrows $y'_{i}\colon {Y\to \cycle_{n-1}E}$ that satisfy $\cycle_{n-1}p\comp y'_{i}= (\del_{j})_{j}\comp \pr_{i}\comp x\comp q$. Again using the induction hypothesis, we get a regular epimorphism $q_{1}\colon {Y_{1}\to Y}$ and arrows $z_{i}\colon {Y_{1}\to E_{n}}$ satisfying $p_{n}\comp z_{i}=\pr_{i}\comp x\comp q\comp q_{1}$, for every $i\in [n+1]$. By the careful choice of the arrows $y'_{i}$ (using the set $K$), these $z_{i}$ induce the needed morphism $z\colon {Y_{1}\to \cycle_{n}E}$.

Considering the diagram
\[
\xymatrix{0 \ar[r] & {K[p]_{n+1}} \ar@{{ >}->}[rr] \ar[d]_-{(\del_{i})_{i\in [n+1]}} && E_{n+1} \ar[d]_-{(\del_{i})_{i\in [n+1]}} \ar@{-{ >>}}[rr]^-{p_{n+1}} && B_{n+1} \ar[d]^-{(\del_{i})_{i\in [n+1]}} \ar[r] & 0\\
0 \ar[r] & {\cycle_{n}K[p]} \ar@{{ >}->}[rr] && \cycle_{n}E \ar@{-{ >>}}[rr]_-{\cycle_{n}p} && {\cycle_{n}B} \ar[r] & 0,}
\]
we see that~\ref{Diagram-Acyclic-Fibration} is a regular pushout if and only if the morphism 
\[
(\del_{i})_{i}\colon {K[p]_{n+1}\to \cycle_{n}K[p]}
\]
is regular epic. The result now follows from Proposition~\ref{Proposition-Characterization-Acyclic-Object}.
\end{proof}

\begin{proposition}\label{Proposition-Characterization-Acyclic-Fibration}
A simplicial morphism $p\colon E\to B$ is an acyclic fibration if and only if $p$ is a regular epimorphism and for every $n\in \N$, Diagram~\ref{Diagram-Acyclic-Fibration} is a regular pushout.
\end{proposition}
\begin{proof}
This is an immediate consequence of Proposition~\ref{Proposition-Characterization-Reg-Epi-Homology-Iso}: By Proposition~4.4 in~\cite{EverVdL2}, a regular epimorphism is always a Kan fibration; conversely, it easily follows from the definitions that a Kan fibration $p$ such that $H_{0}p$ is a regular epimorphism is also regular epic (cf.~\cite[\S II.3, Proposition 1 ff.]{Quillen}).
\end{proof}

\section{Weak equivalences are homology isomorphisms}\label{Section-Weak-Equivalences-Homology-Isos}
We show that for a semi-abelian category $\Ac$ with enough projectives, the weak equivalences in $\simpA$ are exactly the homology isomorphisms.

Recall from~\cite[Section II.4]{Quillen} that a weak equivalence in the sense of Quillen's Theorem is a simplicial morphism $f\colon A\to B$ such that $\hom (P,f)$ is a weak equivalence in $\Sc \Set$ for every projective object $P$ in $\Ac$; and this is the case when $\hom (P,f)$ induces isomorphisms of the \emph{homotopy groups} 
\[
\pi_{n} \hom (P,f)\colon \pi_{n} (\hom (P,A),x)\to \pi_{n} (\hom (P,B), f\comp x)
\]
for all $n\in \N$ and $x\in \hom (P,A)_{0}=\hom (P,A_{0})$. These homotopy groups are defined in the following manner.

\begin{definition}\label{Definition-Homotopy-Groups-Kan}
Let $K$ be a Kan simplicial set, $x\in K_{0}$ and $n\in \N$. Write $Z_{0} (K,x)=K_{0}$; if $n\geq 1$, denote by $x$ the element $(\sigma_{0}\comp \cdots \comp \sigma_{0}) (x)$ of $K_{n-1}$ and by $Z_{n} (K,x)$ the set $\{z\in K_{n}\,|\, \text{$\del_{i} (z)= x$ for $i\in [n]$}\}$. Two elements $z$ and $z'$ of $Z_{n} (K,x)$ are \emph{homotopic}, notation $z\sim z'$, if there is a $y\in K_{n+1}$ such that
\[
\del_{i} (y) = \begin{cases}
x,  & \text{if $i < n$;}\\
z,& \text{if $i = n$;}\\
z',& \text{if $i = n+1$}.
\end{cases}
\]
The relation $\sim$ is an equivalence relation on $Z_{n} (K,x)$; the quotient $Z_{n} (K,x) /_{\sim}$ is denoted by $\pi_{n} (K,x)$. For $n\geq 1$ this $\pi_{n} (K,x)$ is a group.
\end{definition}

\begin{lemma}\label{Lemma-Forks}\cite[Prop.~3.9]{EverVdL2}, \cite[Corollary 6]{Bourn2001}
In a semi-abelian category, let
\[
\xymatrix{A \ar@<1ex>[r]^-{\del_{1}} \ar@<-1ex>[r]_-{\del_{0}} & B \ar[l]|-{\sigma} \ar[r]^-{e} & C}
\]
be an augmented reflexive graph. Then the following are equivalent:
\begin{enumerate}
\item $e$ is a coequalizer of $\del_{0}$ and $\del_{1}$;
\item $e$ is a cokernel of $\del_{1}\comp \ker \del_{0}$.\noproof 
\end{enumerate}
\end{lemma}

The next lemma is essentially a higher-degree version of Proposition~\ref{Proposition-Characterization-H_{0}}: It shows how a homology object $H_{n}A$ may be computed as a coequalizer of (a restriction of) the last two boundary operators $\del_{n}$ and $\del_{n+1}$.

\begin{lemma}\label{Lemma-Kernel-Pair-of-q_{n}}
Write $S_{n+1}A=\bigcap_{i\in [n-1]}K[\del_{i}\colon A_{n+1}\to A_{n}]\subset A_{n+1}$ and $S_{1}A=A_{1}$, and for $n\in \N$, let $\del_{n}^{-1}Z_{n}A$ denote the inverse image
\[
\xymatrix{{\del_{n}^{-1}Z_{n}A} \ar@{}[rd]|<<{\pullback} \ar@<-1ex>[r]_-{\overline{\del}_{n}} \ar@<1ex>[r]^-{\overline{\del}_{n+1}} \ar[d]_-{\del_{n}^{-1}\bigcap_{i\in [n]} \ker \del_{i}} & Z_{n}A \ar[l]|-{\overline{\sigma}_{n}} \ar[d]^-{\bigcap_{i\in [n]} \ker \del_{i}}\\
S_{n+1}A\ar[r]_-{\del_{n}} & A_{n}}
\]
 of $Z_{n}A=\bigcap_{i\in [n]}K[\del_{i}]\subset A_{n}$ along $\del_{n}\colon {S_{n+1}A\to A_{n}}$ and $\overline{\del}_{n+1}$, $\overline{\sigma}_{n}$ the induced factorization of $\del_{n+1}\colon {A_{n+1}\to A_{n}}$ and $\sigma_{n}\colon {A_{n}\to A_{n+1}}$. Then
\[
\xymatrix{{\del_{n}^{-1}Z_{n}A} \ar@<-.5 ex>[r]_-{\overline{\del}_{n}} \ar@<.5 ex>[r]^-{\overline{\del}_{n+1}} & Z_{n} A \ar@{-{ >>}}[r]^-{q_{n}} & H_{n}A} 
\]
is a coequalizer diagram.
\end{lemma}
\begin{proof}
The composite 
\[
\xymatrix{K[\overline{\del}_{n}] \ar@{{ >}->}[r]^-{\ker \overline{\del}_{n}} & \del_{n}^{-1}Z_{n}A \ar[r]^-{\overline{\del}_{n+1}} & Z_{n}A }
\]
is $d'_{n+1}:N_{n+1}A\to Z_{n}A$; the result now follows from Lemma~\ref{Lemma-Forks}.
\end{proof}

\begin{proposition}\label{Proposition-Homology-and-Homotopy-Groups}
Let $\Ac$ be semi-abelian category. For any $n\geq 1$ and for any projective object $P$ in $\Ac$, there is a bijection $\pi_{n} (\hom (P,A),0)\cong \hom (P, H_{n}A)$, natural in $A\in \Ac$.
\end{proposition}
\begin{proof}
An isomorphism $\pi_{n} (\hom (P,A),0)\to  \hom (P, H_{n}A)$ may be defined as follows: The class with respect to $\sim$ of a map $f\colon {P\to Z_{n}A}$ is sent to $q_{n}\comp f\colon P\to H_{n}A$. This function is well-defined and injective, because by Lemma~\ref{Lemma-Kernel-Pair-of-q_{n}}, $f\sim g$ is equivalent to $q_{n}\comp f=q_{n}\comp g$; it is surjective, since $q_{n}$ is a regular epimorphism and $P$ a projective object.  
\end{proof}

\begin{theorem}\label{Theorem-Weak-Equivalences-are-Homology-Isos}
If $\Ac$ is semi-abelian with enough projectives, weak equivalences in $\Sc\Ac$ and homology isomorphisms coincide. 
\end{theorem}
\begin{proof}
If $f$ is a weak equivalence then all $\pi_{n} \hom (P,f)$ are isomorphisms; but then the morphisms $\hom (P,H_{n}f)$, $n\geq 1$ are iso by Proposition~\ref{Proposition-Homology-and-Homotopy-Groups}. Now the class of functors $\hom (P,\cdot)\colon {\Ac \to \Set}$, $P$ projective, is jointly conservative: It reflects regular epimorphisms by definition of a regular projective object, and it reflects monomorphisms because enough regular projectives exist. Hence $H_{n}f$ is an isomorphism.

Conversely, suppose that $f\colon {A\to B}$ is a homology isomorphism. We may factor it as $f=p\comp i$, a trivial cofibration $i\colon {A\to E}$ followed by a fibration $p\colon {E\to B}$. We just proved that $i$ is a homology isomorphism; as a consequence, so is $p$. Thus we reduced the problem to showing that when a fibration is a homology isomorphism, it is a weak equivalence.

Consider $n\in \N$, $P$ projective in $\Ac$, and $x\colon {P\to E_{0}}$. Then 
\[
p_{n}\comp (\cdot)\colon {\pi_{n} (\hom (P, E),x)\to \pi_{n} (\hom (P,B),p_{0}\comp x)}
\]
is an isomorphism: It is a surjection, because by Proposition~\ref{Proposition-Characterization-Acyclic-Fibration}, any map $b\colon {P\to B_{n}}$ satisfying $\del_{i}\comp b=p_{n-1}\comp x$ induces a morphism $e\colon {P\to E_{n}}$ such that $\del_{i}\comp e= x$. It is an injection, since for $e$, $e'$ in $Z_{n} (\hom (P,E),x)$ with $p_{n}\comp e\sim p_{n}\comp e'$, again Proposition~\ref{Proposition-Characterization-Acyclic-Fibration} implies that $e\sim e'$. 
\end{proof}

In summary: 

\begin{theorem}\label{Theorem-Simp-Model-Structure}
If $\Ac$ is a semi-abelian category with enough projectives then a model category structure on $\simpA$ exists where the weak equivalences are the homology isomorphisms and the fibrations are the Kan fibrations.\noproof 
\end{theorem}

\section*{Acknowledgement}

Thanks to Tomas Everaert for many useful comments and suggestions.

\small

\providecommand{\bysame}{\leavevmode\hbox to3em{\hrulefill}\thinspace}
\providecommand{\MR}{\relax\ifhmode\unskip\space\fi MR }
\providecommand{\MRhref}[2]{%
  \href{http://www.ams.org/mathscinet-getitem?mr=#1}{#2}
}
\providecommand{\href}[2]{#2}

\end{document}